\documentclass[letterpaper,11pt]{amsart}
\usepackage{latexsym,array,delarray,amsthm,amssymb,epsfig}%eufrak

\theoremstyle{plain}
\newtheorem{thm}{Theorem}[section]
\newtheorem{lemma}[thm]{Lemma}
\newtheorem{prop}[thm]{Proposition}
\newtheorem{cor}[thm]{Corollary}

\theoremstyle{definition}
\newtheorem{defn}[thm]{Definition}
\newtheorem{ex}[thm]{Example}

\newtheorem{alg}[thm]{Algorithm}

\theoremstyle{remark}
\newtheorem*{rmk}{Remark}
\newtheorem*{obs}{Observation}

\newcommand{\zz}{\mathbb{Z}}
\newcommand{\nn}{\mathbb{N}}
\newcommand{\pp}{\mathbb{P}}

\newcommand{\kk}{\mathbb{K}}

\newcommand{\bfi}{\mathbf{i}}

\newcommand{\bft}{\mathbf{t}}

\newcommand{\bfv}{\mathbf{v}}

\newcommand{\bfx}{\mathbf{x}}
\newcommand{\bfy}{\mathbf{y}}

\newcommand{\ass}{\mathrm{Ass}}

\newcommand{\ind}{\mbox{$\perp \kern-5.5pt \perp$}}

\newcommand{\Sec}{{\rm Sec}}
\newcommand{\al}{\ensuremath{\alpha}}
\newcommand{\bal}{\ensuremath{\boldsymbol{\alpha}}}
\newcommand{\be}{\ensuremath{\beta}}
\newcommand{\bbe}{\boldsymbol{\beta}}

\newcommand{\bga}{\ensuremath{\boldsymbol{\gamma}}}
\newcommand{\bx}{\ensuremath{{\bf x}}}

\begin{document}

\title{Prolongations and computational algebra}

\author{Jessica Sidman}
\address{Department of Mathematics and Statistics, Mount Holyoke
  College, South Hadley, MA 01075}
\email{jsidman@mtholyoke.edu}
\author{Seth Sullivant}
\address{Department of Mathematics and Society of Fellows, Harvard University, Cambridge, MA 02138}
\email{seths@math.harvard.edu}
\thanks{}

\begin{abstract}
We explore the geometric notion of prolongations in the setting of computational algebra, extending results of Landsberg and Manivel which relate prolongations to equations for secant varieties.  We also develop methods for computing prolongations which are combinatorial in nature.  As an application, we use prolongations to derive a new family of secant equations for the binary symmetric model in phylogenetics.
\end{abstract}

\maketitle

%%%%%%%%%%%%%%%%%%%%%%%%%%%%%%%%%%%%%
%%%%%%%%%%%%%%%%%%%%%%%%%%%%%%%%%%%%%
%%%%%%%%%%%%%%%%%%%%%%%%%%%%%%%%%%%%%
%%%%%%%%%%%%%%%%%%%%%%%%%%%%%%%%%%%%%

\section{Introduction}

The notion of \emph{prolongation} originated with Cartan in the context of differential geometry \cite{BCGGG, G, IL}.   We give the simplest formulation of the definition below.  Since we are using differential operators, we will assume that our ground field $\kk$ has characteristic zero.

\begin{defn}\label{defn: differential def}
Let $S = \kk[x_1, \ldots, x_n],$ and let $A \subseteq S_d$ be a vector space of polynomial forms of degree $d.$  The $r$-th \emph{prolongation} of $A,$ denoted by $A^{(r)},$ is
\[
\left\{f \in S_{d+r} \mid \frac {\partial^r f}{\partial \bfx^{\bbe}} \in A \ \mathrm{for  \ all} \ \bbe \in \nn^n \ \mathrm{with} \ |\bbe| = r \right\}.
\]
\end{defn}
In this paper our interest lies not in the role of prolongation in differential geometry, but instead in exploring the applications of purely algebraic reformulations of the definition to three areas: algebraic geometry, commutative algebra,  and phylogenetics.   In particular, we will explain and generalize results of Landsberg and Manivel \cite{LM1} connecting prolongations and secant varieties, as well as apply these ideas to the computation of some nontrivial secant equations arising in phylogenetics.

 Recall that if $X\subset \pp^{n-1}$ is a projective variety, $\Sec^r(X) = X^{\{r\}}$ is the Zariski closure of the union of all  $r-1$ planes spanned by $r$ points in $X.$  Let $I = I(X)$ and suppose that the smallest degree of a minimal generator of $I$ is $d$.   
If $A = I_d,$ the significance of the prolongation $A^{(r)}$ comes from connections to secant ideals of $I.$

\begin{thm}\label{thm: sec}
Let  $X \subseteq \pp^{n-1}$ be a variety over an algebraically closed field of characteristic $0$, with $I = I(X) \subseteq \left<x_1, \ldots, x_n \right>^d$, and $A = I_d$.  Then $A^{((d-1)(r-1))}$ is the degree $r(d-1)+1$ piece of the ideal of the secant variety $X^{\{r\}}$. 
\end{thm}

Theorem \ref{thm: sec} generalizes Lemma 2.2 of \cite{LM1} which concerned prolongations of spaces of quadratic forms.  In the subsequent paper \cite{LM2}, Landsberg and Manivel allude to a generalization of their lemma for higher degrees, but never give a precise statement. At the heart of the proof of Theorem \ref{thm: sec} are connections relating prolongations to polarizations of homogeneous forms and to symbolic powers of ideals. Relationships between symbolic powers and secant ideals are not new; for example containment of secant ideals in symbolic powers appears explicitly in Proposition 2.1 of Catalano-Johnson in \cite{cj} and Corollary 4.6 of Vermeire in \cite{v} shows that if $X$ is smooth, a graded piece of a high enough symbolic power cuts out its secant line variety set-theoretically.

One reason that Theorem \ref{thm: sec} is useful is that it allows for the straightforward computation of equations that belong to the secant ideals $I(X^{\{r\}})$ using linear algebra.  Note that by \cite{Simis2000}, $r(d-1) +1$ is the lowest degree in which there can exist a nonzero polynomial in $I(X^{\{r\}})$.  Even if $I$ is generated by $A = I_d$, it need not be the case that $A^{((d-1)(r-1))}$ generates $I(X^{\{r\}}).$  In spite of this, in several instances of practical interest, the prolongation provides many nontrivial equations in secant ideals that are difficult to derive directly from the definition of the secant variety.  Thus,  prolongation at least brings us one graded piece closer to understanding secant ideals.  

We conclude the paper by returning to our original motivation which was the relationship to algebraic statistics, where algebraic varieties are interpreted as statistical models.  Passing to the secant variety in algebraic geometry amounts to taking a mixture model in statistics.   We will use prolongations as a tool for describing non-trivial secant equations for the binary symmetric model in phylogenetics, which has received attention recently in the algebraic geometry community \cite{Buczynska2006, Sturmfels2005}.

This paper is organized as follows.  We describe several equivalent defintions of prolongation in Section 2, ending with a proof of the relationship between prolongations and symbolic powers.   In Section \ref{sec:comp}, we describe algorithms for computing prolongations.  Theorem \ref{thm: sec}, which connects prolongations and secant equations will follow from results in  Section \ref{sec:secant}.  Part of our proof follows along the lines of ideas from \cite{LM1} and \cite{LM2}, which we attempt to make more explicit, and part of the proof depends on a new application of the join of two ideals.
We derive some nontrivial secant equations for the binary symmetric models in Section \ref{sec:phylo}.

\subsection*{Acknowledgements}
We would like to thank Bernd Sturmfels for suggesting the problem of working on prolongations to us and encouraging our collaboration.  We also thank Aldo Conca, David Cox, J.M. Landsberg, Jason Morton, and Peter Vermeire for helpful conversations and Mike Stillman for help with \emph{Macaulay 2} \cite{GS} which played an important role in our understanding of examples.  Experimentation with MAGMA \cite{magma} and Maple was also vital to this work.  Sidman was partially supported by NSF grant DMS-0600471 and the Clare Boothe Luce Program.  Finally, we are grateful to the Fields Institute for providing a wonderful working atmosphere at an important stage in the project.

%%%%%%%%%%%%%%%%%%%%%%%%%%%%%%%%%%%%%
%%%%%%%%%%%%%%%%%%%%%%%%%%%%%%%%%%%%%
%%%%%%%%%%%%%%%%%%%%%%%%%%%%%%%%%%%%%
%%%%%%%%%%%%%%%%%%%%%%%%%%%%%%%%%%%%%

\section*{Notations and Conventions} \label{sec:conventions}

Since we are working with differential operators, unless explicitly stated, we assume that $\kk$ is a field of characteristic zero.  In situations where we need $\kk$ to be algebraically closed, we explicitly state this.

All varieties $X$ are projective and reduced, but they need not be irreducible; that is, we assume that the ideal $I(X)$ is homogeneous and radical, but not necessarily prime.

Let $\nn$ denote the non-negative integers.  If $\boldsymbol{\alpha} = (\alpha_1, \ldots, \alpha_n) \in \nn^n,$  then $\boldsymbol{\alpha}! = \alpha_1! \cdots \alpha_n!.$  If $| \boldsymbol{\alpha} | = d,$ we let ${d \choose \boldsymbol{\alpha}}$ denote the multinomial coefficient with parts $\alpha_1, \ldots, \alpha_n.$  We write $\bf{1}$ to denote a vector in which each coordinate is 1.  A monomial in the polynomial ring $\kk[x_1, \ldots, x_n] = \kk[\bx]$ is given by an element $\bal \in \nn^n$ where $\bx^{\bal} = x_1^{\al_1}\cdots x_n^{\al_n}.$

%%%%%%%%%%%%%%%%%%%%%%%%%%%%%%%%%%%%%
%%%%%%%%%%%%%%%%%%%%%%%%%%%%%%%%%%%%%
%%%%%%%%%%%%%%%%%%%%%%%%%%%%%%%%%%%%%
%%%%%%%%%%%%%%%%%%%%%%%%%%%%%%%%%%%%%

\section{Prolongations:  Equivalent Definitions}\label{sec:equiv}

In this section we focus on making the translations between purely algebraic descriptions of prolongation and Definition \ref{defn: differential def} explicit.  The importance of the algebraic definitions is that they allow us to connect the notion of prolongations to commutative algebra and algebraic geometry, specifically, to symbolic powers and equations defining secant varieties.

In the interest of keeping our introduction to prolongations self-contained, we briefly review the definition of the symmetric algebra in terms of the tensor algebra as well as the basics of polarization although this material will be well-known to some.  We define prolongations from the quotient algebra point of view in Section \ref{sec: sym}.  We translate the definition into the language of polarization in Section \ref{sec: polar}.  Finally, we illustrate how thinking of prolongation in terms of polarization leads to connections with symbolic powers.

\subsection{Prolongations and the symmetric algebra}\label{sec: sym}
In this section we will recast the definition of prolongation in terms of the symmetric algebra, viewed as a quotient of the tensor algebra. 

 Let $V$ be a finite dimensional vector space over a field $\kk,$ and $V^*$ denote its dual.  We follow the conventions of Appendix B of \cite{FH}.  The reader may also want to consult Appendix A2 of \cite{E} or Chapter 1 of \cite{W}.

We let $T = \bigoplus_{d \geq 0} T^dV^*$ be the tensor algebra of $V^*,$ where $T^dV^*$ denotes the tensor product of $V^*$ with itself $d$ times.  The symmetric algebra $S,$ on $V^*$ is defined to be the quotient of $T$ by the ideal $\langle x \otimes y - y \otimes x \mid x, y \in V^* \rangle.$  If we pick a basis $\bx = (x_1, \ldots, x_n),$ for $V^*$, then we may identify $S$ with $\kk[x_1, \ldots, x_n],$ the homogeneous coordinate ring of $\pp V.$  A monomial $\bx^{\bal}$ of degree $d$ represents the equivalence class of tensors which map to it under the canonical projection $T \to S.$

The \emph{co-multiplication} or \emph{diagonal} map $S^{d+r}V^* \to S^dV^* \otimes S^rV^*$ will be important in what follows.
First, we describe it using the intrinsic point of view of Appendix A2.4 in \cite{E}.  Recall that the diagonal map $\Delta:S \to S \otimes S$ sends  $x \in V^*\mapsto x\otimes 1+ 1 \otimes x.$ We get a map $\Delta_{d,r}: S^{d+r}V^* \to S^dV^* \otimes S^rV^*$ by restricting the diagonal map to $S^{d+r}V^*$ and composing this with the projection to $S^dV^* \otimes S^rV^*.$ For example, since
\begin{align*}
\Delta(x^2y)& = (x \otimes 1+ 1 \otimes x)^2(y \otimes 1+1 \otimes y) \\
 & =  x^2y \otimes 1 +2xy \otimes x + y\otimes x^2 + x^2 \otimes y + 2x \otimes xy +1 \otimes x^2y,
\end{align*}
we see that the co-multiplication map $\Delta_{2,1}$ sends $x^2y$ to $2xy \otimes x+x^2 \otimes y \in S^2V^* \otimes S^1V^*.$ 
 
Following pg.~5 of \cite{W} we can also describe the co-multiplication map by its action on monomials in terms of our basis $\bx.$  If $i_1 \leq \cdots \leq i_{d+r},$ then
\[
x_{i_1}\cdots x_{i_{d+r}} \mapsto \sum x_{i_{\sigma(1)}} \cdots x_{i_{\sigma(d)}}\otimes x_{i_{\sigma(d+1)}}\cdots  x_{i_{\sigma(d+r)}}
\]
where we sum over all  permutations $\sigma$ of $d+r$ elements such that $\sigma(1) < \cdots < \sigma(d)$ and $\sigma(d+1) < \cdots < \sigma(d+r).$ 

The reason that co-multiplication appears in connection with prolongation is that it is closely related to partial differentiation.

\begin{lemma}
If $F \in S^{d+r}V^*,$ then $\Delta_{d,r}(F)=\displaystyle \sum_{|\bbe| = r} \frac{1}{\bbe !}\frac{\partial^r F}{\partial \bx^{\bbe}} \otimes \bx^{\bbe}.$
\end{lemma}
\begin{proof}
By linearity it suffices to assume that $F$ is a monomial $F = \bfx^{\bal}.$ The projection of $\Delta(\bfx^{\bal})$ to $S^dV^* \otimes S^rV^*$ is the sum of all of the terms of $\Delta(\bfx^{\bal})$ which can be written in the form $- \otimes \bfx^{\bbe}$ with $|\bbe| = r.$  Since 
\[
\Delta(\bfx^{\bal}) = \prod_{i=1}^n (x_i \otimes 1 + 1 \otimes x_i)^{\alpha_i},\]
there will be ${\alpha_1 \choose \beta_1} \cdots {\alpha_n \choose \beta_n}$ terms in the product of the form  $- \otimes \bfx^{\bbe},$ all equal to $\bfx^{\bal -\bbe} \otimes \bfx^{\bbe}.$  But
\[
{\alpha_1 \choose \beta_1} \cdots {\alpha_n \choose \beta_n} \bfx^{\bal -\bbe} \otimes \bfx^{\bbe} = \frac{1}{\bbe !}\frac{\partial^r \bfx^{\bal}}{\partial \bx^{\bbe}} \otimes \bx^{\bbe}
\] 
\end{proof}
 We can use co-multiplication to see that the following algebraic definition of prolongation given in Section 2.1.3 of \cite{LM1} is equivalent to Definition \ref{defn: differential def}.

\begin{lemma}\label{lem: prolong}
If $A \subset S^dV^*,$ then 
$A^{(r)} =  \left(A \otimes S^{r}V^*\right) \cap S^{d+ r}V^*.$
\end{lemma}
\begin{proof}
Note that $\Delta_{d,r}$ maps $F \in S^{d+r}V^*$ to an element of the form
\[
 \sum_{|\bbe|=r} F_{\bbe} \otimes {\bf x}^{\be} \in S^dV^* \otimes S^rV^*.
\] This is in $A \otimes S^rV^*$ if and only if each $F_{\bbe} \in A,$ and  by the previous lemma, $F_{\bbe
} = \frac1{\bbe !}\frac {\partial^r F}{\partial \bfx^{\bbe}}.$ 
\end{proof}

%%%%%%%%%%%%%%%%%%%%%%%%%%%%%%
%%%%%%%%%%%%%%%%%%%%%%%%%%%%%%
%%%%%%%%%%%%%%%%%%%%%%%%%%%%%%
%%%%%%%%%%%%%%%%%%%%%%%%%%%%%%

\subsection{Prolongations and polarization}\label{sec: polar}

In this section we explain the connection between prolongation and \emph{polarization}.  Polarization, which arose in classical invariant theory \cite{Weyl}, is the higher degree analog of associating a symmetric bilinear form to a quadratic form and is closely related to the representation of a homogeneous form as an element of the tensor algebra.  The notion of polarization is also used in \S 1 of Chapter VI in \cite{ACGH} in connection with secant varieties of curves.  What is significant for us is that thinking of prolongation in terms of polarization opens the door to connections with symbolic powers and with secant varieties.

\begin{defn}
Suppose that $F$ is a homogeneous polynomial of degree $d$ in $\kk[\bx]$ where $\bx = (x_1, \ldots, x_n).$  For each $i =1, \ldots, d$ we introduce a new set of $n$ variables $\bx_i = (x_{i1}, \ldots, x_{in}).$  We also introduce an auxiliary set of variables $\bft =(t_1, \ldots, t_d).$ The \emph{polarization} of $F,$ denoted ${\bf F}(\bx_1, \ldots, \bx_d),$ is the coefficient of $\bft^{\bf 1}$ in the expansion of $F(t_1\bx_1+ \cdots +t_d\bx_d)$ as a polynomial in $\bft$.
\end{defn}

\begin{ex}
Let $F(\bx) = x_1^2x_2.$  We compute
\begin{align*}
F(t_1\bx_1+t_2\bx_2+t_3\bx_3) & = (t_1 x_{11}+t_2 x_{21} +t_3 x_{31})^2(t_1 x_{12} + t_2 x_{22} + t_3 x_{32})\\
&= t_1^3x_{11}^2x_{12}+t_2^3x_{21}^2x_{22} + \cdots\\
& + 3t_1^2t_2(x_{11}^2x_{22}+x_{11}x_{21}x_{12})+ \cdots \\
& + t_1t_2t_3(2x_{11}x_{21}x_{32} + 2x_{11}x_{31}x_{22} + 2x_{21}x_{31}x_{12}).
\end{align*}
We see that ${\bf F}(\bx_1, \bx_2, \bx_3) = 2x_{11}x_{21}x_{32} + 2x_{11}x_{31}x_{22} + 2x_{21}x_{31}x_{12}.$
\end{ex}
\begin{lemma}
If $F(\bx)$ is a homogeneous form of degree $d,$ then
\begin{enumerate}
\item  \[F(t_1\bx_1+ \cdots +t_d\bx_d) = \sum_{|\bbe| = d} \frac{\bft^{\bbe}}{\bbe !} {\bf F}(\bx_1^{  \beta_1 }, \ldots, \bx_d^{\beta_d})\]
where $\bx_i^{\beta_i}$ means that the set of variables $\bx_i$ is repeated $\beta_i$ times.

\item  ${\bf F}(\bx_1, \ldots, \bx_d)$ is linear in each set of variables $\bfx_i.$

\item   ${\bf F}(\bx_1, \ldots, \bx_d)$ is symmetric in the $\bfx_i.$ (If $\sigma$ is a permutation of $d$ elements, then  $ {\bf F}(\bx_1, \ldots, \bx_d)={\bf F}(\bx_{\sigma(1)}, \ldots, \bx_{\sigma(d)}).$)

 \item ${\bf F}(\bx, \ldots, \bx) = d!F(\bx)$
\end{enumerate}
\end{lemma}

\begin{proof}
Note that it is enough to prove the stated claims in the case where $F(\bx) = \bx^{\bal}.$  In this case
\begin{equation}\label{eq:  prod}
F(t_1\bx_1+ \cdots +t_d\bx_d) = \prod_{j=1}^n (t_1 x_{1j} + \cdots + t_d x_{dj})^{\alpha_j}.
\end{equation}

For part (1), recall that ${\bf F}(\bfx_1, \ldots, \bfx_d)$ is the coefficient of $\bft^{\bf 1}$ in the product consisting of $d$ factors of the form $(t_1x_{1j}+ \cdots + t_d x_{dj})$ as above.  From this definition, we see that the coefficient of $\bft^{\bf 1}$ is a sum of $d!$ monomials (counted with multiplicity) which correspond to the $d!$ ways of choosing one term per factor, where each $t_j$ is chosen exactly once.

We can construct the $d!$ monomials which are coefficients of $\bft^{\bf 1}$ as follows.  Let the $d$ factors in Equation (\ref{eq:  prod}) be $F_1, \ldots,  F_d.$  Fix a subset $I \subset [d]$ of size $\beta.$  A monomial coefficient of $\bft^{\bf 1}$ is gotten by choosing $\beta$ terms of the form $t_kx_{kj}$ from among the factors $F_s$ with $s \in I$ where $k\in I$ and the $k$ are all distinct, and $d-\beta$ terms $t_{\ell}x_{\ell j}$  from among the factors $F_t$ with $t \notin I$ where $\ell \notin I$ and the $\ell$ are all distinct.  Since there are ${d \choose \beta}$ ways to choose the set $I,$ $\beta!$ ways to make our construction of a monomial from the factors $F_s$ with $s \in I$ and $(d-\beta)!$ ways to construct a monomial from the factors $F_t$ with $t \notin I,$ we see that we get each monomial coefficient of $\bft^{\bf 1}$ in this way.

Let us consider ${\bf F}(\bx_1^{  \beta_1 }, \ldots, \bx_d^{\beta_d})$ where $\bfx_i$ is repeated $\beta_i$ times.  We can pass from  ${\bf F}(\bfx_1, \ldots, \bfx_d)$ to ${\bf F}(\bx_1^{  \beta_1 }, \ldots, \bx_d^{\beta_d})$  by computing the result of repeating the $i$-th set of  variables $\beta_i$ times successively for each $i.$

Assume that $I \subset [d]$ is of size $\beta_i$, and that for each $k \in I,$ $t_kx_{kj}$ is replaced by $t_ix_{ij}.$  In ${\bf F}(\bfx_1, \ldots, \bfx_d)$ there are $\beta_i!$ ways of choosing terms of the form $t_kx_{kj}$ from the factors $F_s$ with $s \in I$ with $k\in I$  and the $t_k$ all distinct, but if we replace each  $t_kx_{kj}$ by $t_ix_{ij},$ all of these $\beta_i!$ choices look the same.

Now let us turn to the consideration of the coefficent of $\bft^{\bbe}$ in Equation (\ref{eq:  prod}).  Note that there is only one way of choosing $\beta_i$ terms of the form $t_ix_{ij}$ from among the factors $F_k$ with $k \in I.$  Considering each of the $d$ sets of variables in turn, we see that  if we substitute the $i$th set of variables $\beta_i$ times in ${\bf F}(\bfx_1, \ldots, \bfx_d)$ we will see each monomial that appears as a coefficient of $\bft^{\bbe}$ repeated $\bbe!$ additional times.  Therefore, the coefficient of $\bft^{\bbe}$ is $\frac1{\bbe!} {\bf F}(\bx_1^{  \beta_1 }, \ldots, \bx_d^{\beta_d}).$

Parts (2) and (3) follow immediately from the definition of ${\bf F}(\bx_1, \ldots, \bx_d)$ as the coefficient of $\bft^{\bf 1}$ in Equation (\ref{eq:  prod}).  Part (4) follows by computing ${\bf F}(\bx, \ldots, \bx)$ as the coefficient of $\bft^{\bf 1}$ in 
\[F(t_1 \bx+ \cdots + t_d \bx) = \prod_{i=1}^n ((t_1+\cdots +t_d)x_i)^{\alpha_i} = (t_1 + \cdots + t_d)^d\bx^{\bal}.\]
\end{proof}

\begin{obs}
Recall that we may identify the elements of $S^dV^*$ with elements of $T^dV^*$ that are invariant under the action of the symmetric group on $d$ letters.  This point of view is especially important in \cite{LM1} and \cite{LM2}.  Explicitly, if $v_i \in V^*,$ then 
\[v_1 \cdots v_d \mapsto \sum_{\sigma \in S_d} v_{\sigma(1)} \otimes \cdots \otimes v_{\sigma(d)}.\]  Note that the image of a monomial $m \in S^dV^*$ in $T^dV^*$ is a weighted sum over all of the monomials in the coset represented by $m,$ and $m$ is $1/d!$  times the projection of this element of $T^dV^*$  into $S^dV^*.$   Therefore, we can easily pass from the expression of an element of $S^dV^*$ as a $d$-tensor that is invariant under the action of the symmetric group on $d$ letters to its polarization; we just write the elements appearing in the $i$-th factor in the tensor in terms of $\bfx_i$ and erase the tensor symbols.  For example,

\[
x_1^2x_2 \mapsto 2(x_1 \otimes x_1 \otimes x_2 + x_1 \otimes x_2 \otimes x_1 +x_2 \otimes x_1 \otimes x_1) \mapsto 2(x_{11}x_{21}x_{32} +x_{11}x_{22}x_{31}+x_{12}x_{21}x_{31}).
\]
\end{obs}

The following elementary lemma describes relationships between polarization and partial differentiation.

\begin{lemma}\label{lem: polar}
Let $F$ be a homogeneous polynomial of degree $d+r.$
\begin{enumerate}
\item  If $F = {\bf x}^{\bal},$ then 
\[ 
{\bf F}({\bf x}, \ldots, {\bf x}, {\bf y}, \ldots, {\bf y}) = \boldsymbol{\alpha}!\sum_{ \boldsymbol{\beta} \in \nn^n, |\boldsymbol{\beta}| = r} {d \choose \boldsymbol{\alpha-\beta}} {r \choose \boldsymbol{\beta}}{\bf x}^{\boldsymbol{\alpha-\beta}}{\bf y}^{\boldsymbol{\beta}}.
\] 

\item  If $\boldsymbol{\beta} \in \nn^n$ with $|\boldsymbol{\beta}| = r,$ then
\[\label{eq: coeff}
d!r!\frac{\partial^r F}{\partial {\bf x}^{\boldsymbol{\beta}}} 
= \frac{\partial^r {\bf F}({\bf x}, \ldots, {\bf x}, {\bf y}, \ldots, {\bf y})}{\partial {\bf y}^{\boldsymbol{\beta}}}.
\]

\end{enumerate}
In both expressions we assume there are $d$ copies of $\bfx$ and $r$ copies of $\bfy$.
\end{lemma}

\begin{proof}
\noindent (1) When we polarize $F$, we get a $(d+r)$-linear symmetric form in $d+r$ sets of variables ${\bf x}_i.$  Each of the ${d+r \choose \boldsymbol{\alpha}}$ distinct monomials in ${\bf F}({\bf x}_1, \ldots, {\bf x}_{d+r})$ appears with coefficient $\boldsymbol{\alpha}!.$  For any choice of $\boldsymbol{\beta},$ there will be ${d \choose \boldsymbol{\alpha-\beta}}{r \choose \boldsymbol{\beta}}$ distinct monomials which will agree (and have ${\bf y}^{\boldsymbol{\beta}}$ as a factor) when the first $d$ sets of variables are all set to ${\bf x}$ and the last $r$ are set to ${\bf y}.$
\bigskip

\noindent (2) It is enough to prove the result for an arbitrary monomial of degree $d+r.$  Assume, without loss of generality, that $F = {\bf x}^{\boldsymbol{\alpha}}.$  Using (1) we see that the coefficient of ${\bf y}^{\boldsymbol{\beta}}$ in ${\bf F}({\bf x}, \ldots, {\bf x}, {\bf y}, \ldots, {\bf y})$ is

\[
\boldsymbol{\alpha}!{d \choose \boldsymbol{\alpha-\beta}} {r \choose \boldsymbol{\beta}}{\bf x}^{\boldsymbol{\alpha-\beta}} = \frac{d!r!\boldsymbol{\alpha}!}{(\boldsymbol{\alpha-\beta})!\boldsymbol{\beta}!}{\bf x}^{\boldsymbol{\alpha-\beta}}.
\]
Therefore, we see that taking partial derivatives with respect to ${\bf y}^{\boldsymbol{\beta}}$ yields

\[d!r!\frac{\boldsymbol{\alpha}!}{(\boldsymbol{\alpha-\beta})!} {\bf x}^{\boldsymbol{\alpha-\beta}} = d!r!\frac{\partial^r F}{\partial {\bf x}^{\boldsymbol{\beta}}}.
\]

\end{proof}
The next lemma is a modification of an observation in \cite{LM1}.  (See also the discussion after Corollary 3.2 in \cite{LM2}.)

\begin{lemma}\label{lem: polar2}
Let $A \subseteq S^dV^*,$ and $F \in S^{d+r}V^*$ be a homogeneous polynomial with polarization ${\bf F}({\bf x}, \ldots, {\bf x}, {\bf y}, \ldots,{\bf y}),$ of degree $d$ in the ${\bf x}$-variables.  The following are equivalent:
\begin{enumerate}
\item  $F$ is in $A^{(r)}.$
\item  Every coefficient of ${\bf F}$ as a polynomial in the ${\bf y}$-variables is in $A.$

\item Every coefficient of 
${\bf F}({\bf x}, \ldots, {\bf x}, {\bf y}_1, \ldots, {\bf y}_r),$ viewed as a polynomial in all of the ${\bf y}$-variables, is in $A.$

\item  ${\bf F}({\bf x}, \ldots, {\bf x}, {\bf v}, \ldots,{\bf v}) \in A$ for every choice of ${\bf v} \in V.$

\end{enumerate}
 
\end{lemma}

\begin{proof}
First we show the equivalence of (1) and (2).  We know that $F \in A^{(r)}$ if and only if $\frac{\partial^r F}{\partial \bfx^{\bbe}} \in A$ for every $\bbe \in \nn^n$ with $|\bbe| = r.$  But $\frac{\partial^r F}{\partial \bfx^{\bbe}}$ is just $\bbe!$ times the coefficient of $\bfy^{\bbe}$ in ${\bf F}$ by part (2) of Lemma \ref{lem: polar}.

The equivalence of (2) and (3) follows because the coefficient of the monomial ${\bf y}^{\boldsymbol{\beta}}$ in the polynomial ${\bf F}({\bf x}, \ldots, {\bf x}, {\bf y}, \ldots, {\bf y}),$ which has degree $r$ in the ${\bf y}$-variables, is ${r \choose \boldsymbol{\beta} }$ times the coefficient of some monomial in the $r$ sets of ${\bf y}$-variables in  ${\bf F}({\bf x}, \ldots, {\bf x}, {\bf y}_1, \ldots, {\bf y}_r).$

For the equivalence of (2) and (4), note that we are working over an infinite field.   If we choose ${d+r+n-1 \choose d+r}$ sufficiently generic points ${\bf v}_i,$ the vectors in indeterminates $F_{\bal}$ of the form $\sum_{|\bal| = r} F_{\bal}{\bf v}_i^{\bal}$
will be linearly independent.  Hence, they are all in $A,$ if and only if every $F_{\bal} \in A.$
\end{proof}

Lemma \ref{lem: polar2} allows us to prove the following result about the prolongations and ideals.

\begin{thm}\label{thm: diff}
Let $A$ be a subspace of $S^dV^*$ and let $I \subset \kk[\bx]$ be the ideal generated by $A.$  Then $F \in S^{d+r}V^*$ is in $A^{(r)}$ if and only if $\frac{\partial^r F}{\partial {\bf x}^{\boldsymbol{\beta}}} \in I$ for every $\boldsymbol{\beta} \in \nn^n$ with $|\boldsymbol{\beta}| =k\leq r.$ 
\end{thm}
\begin{proof}
One inclusion follows immediately from Definition \ref{defn: differential def}.  For the opposite inclusion assume that $F \in A^{(r)}.$   Write  ${\bf F}({\bf x}, \ldots, {\bf x}, {\bf y}_1, \ldots, {\bf y}_r)$ as a polynomial in the variables ${\bf y}_i:$
\[
\sum_{\boldsymbol{\beta} \in \nn^{nr}, |\boldsymbol{\beta}| = k}F_{\boldsymbol{\beta}}({\bf x}, \ldots, {\bf x}) {\bf y}^{\boldsymbol{\beta}}.
\]
The symbol ${\bf y}$ above stands for the vector of vectors ${\bf y} = ({\bf y}_1, \ldots, {\bf y}_{r}).$

Part (3) of Lemma \ref{lem: polar2} tells us that each $F_{\boldsymbol{\beta}}({\bf x}, \ldots, {\bf x})$ is in $A.$  Now set the variables ${\bf y}_1, \ldots, {\bf y}_{r-k}$ equal to ${\bf x}.$  Each 
\[F_{\boldsymbol{\beta}}({\bf x}, \ldots, {\bf x}) {\bf y}^{\boldsymbol{\beta}} \mapsto F_{\boldsymbol{\beta}}({\bf x}, \ldots, {\bf x}) {\bf x}^{\boldsymbol{\alpha}'}{\bf y}^{\boldsymbol{\beta - \alpha}}\]
where $\boldsymbol{\alpha}' \in \nn^n$ has $i$-th coordinate $\sum_j \alpha_{ij}.$  From this, we can see that the coefficients of the monomials in the remaining ${\bf y}_i$-variables are all in $I.$   Setting the remaining ${\bf y}_i$-variables all equal to ${\bf v},$ we see that the coefficients of the monomials in ${\bf v}$ will be in $I$ and may be interpreted as partial derivatives of $F$ of order $k$ (up to a scalar) by part (2) of Lemma \ref{lem: polar}.  We conclude that all partial derivatives of order $k$ are in $I$ for all $k \leq r.$ 
\end{proof}

Recall that if $I \subset S$ is a radical ideal with $\ass(I) = \{P_1, \ldots, P_m\},$ then the $r$-th \emph{symbolic power} of $I$ is defined to be 
\[ I^{(r)} := I^r [W^{-1}] \cap S,\]
where $W = S \backslash (P_1 \cup \cdots P_m).$  We define the $r$-th \emph{differential power} of $I$ to be \[I^{<r>} = \left\{  f \in S \mid \frac{\partial^{|\bbe|} f}{\partial \bx^{\bbe}} \in I  \, \, \,  \mbox{for all} \, \,    |\bbe| \leq r-1 \right\}.\]  If $\kk$ is algebraically closed of characteristic 0 and $I$ is prime, then by the theorem of Zariski and Nagata $I^{(r)} = I^{<r>}.$
See Theorem 3.14 in \cite{E} for a discussion of the proof and pointers to a more general statement in characteristic $p.$  

The theorem of Zariski and Nagata also holds for radical ideals.  While we found this statement in the literature, we could not find its proof, so we include one for completeness.

\begin{cor}
If $I$ is a radical ideal over an algebraically closed field of characteristic zero, then $I^{(r)} = I^{<r>}.$ 
\end{cor}
\begin{proof}
Suppose that $\ass(I) = \{P_1, \dots, P_m\}.$  It is easy to see that $I^{<r>} = \cap P_i^{<r>},$ so by the theorem of Zariski and Nagata it suffices to show that $I^{(r)} = \cap P_i^{(r)}.$  Since $I$ is radical, by prime avoidance,  $\ass(I^{(r)}) = \ass(I),$ and the $P_i$-primary component of $I^{(r)}$ is $(I^{(r)})_{P_i}\cap S.$  But, as localization commutes with products and intersections, we have
\[
(I^{(r)})_{P_i} =( I^r [W^{-1}] \cap S)_{P_i} = (I^r)_{P_i} \cap S_{P_i} = (I_{P_i})^r = ((P_i)_{P_i}))^r = (P_i^r)_{P_i}.
\]
We see that  $(I^{(r)})_{P_i}\cap S =  (P_i^r)_{P_i} \cap S = P_i^{(r)},$ which completes the proof.
\end{proof}

Thus, we have the following corollary to Theorem \ref{thm: diff}.
\begin{cor}\label{cor:symbpower}
Let $\kk$ be an algebraically closed field of characteristic 0.  Let $I$ be a radical ideal with $I \subset \left<x_1, \ldots, x_n \right>^d$ and let $A = I_d$.   Then the $(d+r)$-th graded piece of $I^{(r+1)}$ is $A^{(r)}.$
\end{cor}

%%%%%%%%%%%%%%%%%%%%%%%%%%%%%%%%%%%%%
%%%%%%%%%%%%%%%%%%%%%%%%%%%%%%%%%%%%%
%%%%%%%%%%%%%%%%%%%%%%%%%%%%%%%%%%%%%
%%%%%%%%%%%%%%%%%%%%%%%%%%%%%%%%%%%%%

\section{Computing Prolongations}\label{sec:comp}

In this section we describe algorithms for computing prolongations that use linear algebra and can be implemented in a computer algebra system.  We also discuss how combinatorial tools can be used to speed up the computations by reducing the dimensions of the intermediate vector spaces that need to be computed.  These combinatorial approaches can also be used to determine explicit descriptions of prolongations.

%%%%%%%%%%%%%%%%%%%%%%%%%%%%%%%%%%%%%
%%%%%%%%%%%%%%%%%%%%%%%%%%%%%%%%%%%%%
%%%%%%%%%%%%%%%%%%%%%%%%%%%%%%%%%%%%%
%%%%%%%%%%%%%%%%%%%%%%%%%%%%%%%%%%%%%
\subsection{Algorithms}

We will describe some algorithms for computing prolongations which depend on various implementations of the equivalent definitions from Section \ref{sec:equiv}.  In practice, we will have a basis for $A$ consisting of a set of homogeneous polynomials of degree $d$  and will want to compute $A^{(r)}.$   The crucial step in each of the algorithms for computing prolongations that we describe may be performed by Gaussian elimination.  However, what is easiest to implement depends on the way polynomials are stored since converting polynomials to vectors which may be operated on by the user may be nontrivial in practice in any given computer algebra package. 

\begin{alg}\label{alg: prolong1}~

\noindent  INPUT:  A basis for $A.$

\noindent  OUTPUT:  A basis for $A^{(r)} \subseteq S^{d+r}V^*.$

\noindent STEP 1:  Map a basis for $S^{d+r}V^*$ into $S^dV^* \otimes S^rV^*$ via the co-multiplication map.  

\noindent STEP 2:  Compute the intersection of  $A \otimes S^rV^*,$ with the space constructed in step 1.  Multiplication (just ``erasing'' the tensor symbol) maps a basis for this intersection into $S^{d+r}V^*$ giving a basis for $A^{(r)}.$
\end{alg}

 Unfortunately, the simplest implementation of Algorithm \ref{alg: prolong1} introduces a new set of variables to represent the terms to the right of the tensor symbol, which slows computation. Alternatively, we can exploit the connection between co-multiplication and partial differentiation to avoid introducing a new set of variables.

\begin{alg}~

\noindent INPUT:  A basis for $A.$

\noindent OUTPUT: A basis for $A^{(r)}.$

\noindent STEP 1:  For each $\bbe \in \nn^r,$ with $|\bbe| = r,$ compute $A_{\bbe},$ the space of all forms of degree $d+r$ whose partial derivative with respect to $\bfx^{\bbe}$ is in $A.$

\noindent STEP 2: The intersection of the spaces $A_{\bbe}$ is equal to $A^{(r)}.$  

\end{alg}

 Another alternative would require more extensive programming, but is potentially quite fast.  The complexity of Algorithm \ref{alg: prolong3} is governed by the amount of pre-processing necessary to coordinatize our basis vectors and the Gaussian elimination in STEP 3.

\begin{alg}~ \label{alg: prolong3}

\noindent INPUT:  A basis $B$ for $A.$

\noindent OUTPUT:  A basis for $A^{(r)}.$

\noindent STEP 1:  Fix a term order so that monomials form an ordered basis for $S.$  Abusing notation, let $A$ denote the matrix whose columns are the coefficients of the elements of $B$ with respect to this ordered basis.

\noindent STEP 2:  Let $F$ be the generic form of degree $d+r.$ Let $C$ be a matrix with a column for each $|\bbe|=r.$  The column corresponding to $\bbe$ is the coordinate vector of $\frac{\partial^r F}{\partial \bx^{\bbe}}$ with respect to our ordered basis.  (After a suitable scaling of the basis elements, $C$ is the Catalecticant matrix $C(d,r; n).$)

\noindent STEP 3:  Form the augmented matrix $[A|C].$  The space $A^{(r)}$ is just the space of polynomials $F$ for which the augmented matrix $[A|C]$ is consistent.  We find this space by putting $A$ in reduced-echelon form which will give a linear equation on the entries of $C$ for every zero row to the left of the bar in the augmented matrix.  Solving this system of equations gives the coordinate vectors of a basis for  $A^{(r)}.$
\end{alg}

%%%%%%%%%%%%%%%%%%%%%%%%%%%%%%%%%%%%%
%%%%%%%%%%%%%%%%%%%%%%%%%%%%%%%%%%%%%
%%%%%%%%%%%%%%%%%%%%%%%%%%%%%%%%%%%%%
%%%%%%%%%%%%%%%%%%%%%%%%%%%%%%%%%%%%%
\subsection{Monomial Prolongations}

In this section, we describe the prolongations of vector spaces spanned by monomials.  The monomial case can be solved purely combinatorially and can be used as a tool for reducing the computational burden in the general case.

\begin{prop}
Suppose that $A$ is spanned by monomials.  A monomial $\bfx^{\bal}$ is in $A^{(r)}$ if and only if $\bfx^{\bal- \bbe} \in A$ for all $\bfx^{\bbe}$ dividing $\bfx^{\bal}$ with $|\bbe| = r.$
\end{prop}

\begin{proof}
The differential operator $\frac{\partial^r}{\partial \bfx^{\bbe}}$ maps monomials to monomials.  If $\bfx^{\bbe}$ divides $\bfx^{\bal}$ then $\frac{\partial^r}{\partial \bfx^{\bbe}} \bfx^{\bal} = C \bfx^{\bal - \bbe}$ for a nonzero constant $C$, otherwise $\frac{\partial^r}{\partial \bfx^{\bbe}} \bfx^{\bal}  = 0$.
\end{proof}

An important special case arises when $d = 2$.  In this case, the generators of $A$ have two types:  squarefree pairs $x_ix_j$ and pure powers $x_i^2.$ Let $\sigma \subset [n]$ denote the set of $i$ such that $x_i^2 \in A$.  We define a graph $G(A,r)$ as follows:

\begin{defn}
Let $A$ be a vector space spanned by quadratic monomials from $\kk[\bx].$  For each integer $r > 0,$ we define a graph $G(A,r)$ with $r+2$ vertices for each indeterminate whose square is in $A$, and a single vertex for all other indeterminates. Formally, the vertex set of $G(A,r)$ is the set of all pairs $(i,j)$ with $i \in [n],$ where $j \in [r+2]$ if $i \in \sigma$ and $j = 1$ otherwise.  A pair of vertices $(i_1, j_1)$ $(i_2, j_2)$ is connected by an edge if $x_{i_1}x_{i_2} \in A$.
\end{defn}

The graph $G(A,r)$ can be used to read off the generators of the prolongations $A^{(r)}$.

\begin{cor}\label{cor:blowup}
The induced subgraph of $G(A,r)$ on vertices $(i_1, j_1), \ldots, (i_{r+2}, j_{r+2})$ is a complete graph if and only if $x_{i_1} \cdots x_{i_{r+2}}$ is in the prolongation $A^{(r)}$.
\end{cor}

\begin{proof}
The set of vertices $(i_1, j_1), \ldots, (i_{r+2}, j_{r+2})$ forms a $K_{r+2}$ if and only if $x_{i_k} x_{i_l} \in A$ for all $1 \leq k< l \leq r+2$ if and only if each divisor of $x_{i_1} \cdots x_{i_{r+2}}$ of degree $r$ has quotient in $A$.
\end{proof}

\begin{ex}
Let $A$ be the span of $x_1^2, x_1x_2, x_2x_3, x_2x_4, x_3x_4 \in \kk[x_1, x_2, x_3, x_4].$  To compute $A^{(1)},$ we construct the graph $G(A,1)$:

\[ 
\includegraphics{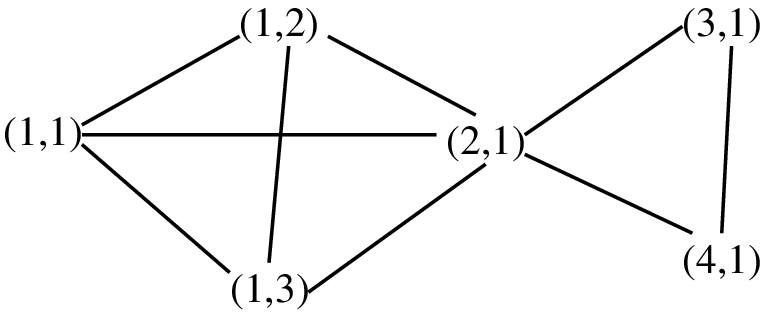}
\]

The graph $G(A,1)$ contains five complete subgraphs $K_3$.  One is the subgraph of $G(A,1)$ on the vertices $(1,1), (1,2), (1,3)$ and this corresponds to the monomial $x_1^3 \in A^{(1)}$.   There are three triangles in the graph of the form $(1,i), (1,j), (2,1)$ and these all correspond to the monomial $x_1^2 x_2 \in A^{(1)}$.  Finally the triangle $(2,1), (3,1), (4,1)$ corresponds to the monomial $x_2x_3x_4 \in A^{(1)}$.  Corollary \ref{cor:blowup} implies that these three monomials span $A^{(1)}$.
\end{ex}

\begin{cor}
Computing prolongations is NP-hard.
\end{cor}

\begin{proof}
Focusing on the case where $M$ is generated by squarefree quadratic monomials, we see that the prolongations are determined by the complete subgraphs in a fixed graph $G(M)$.  In particular, $A^{(r)}$ is empty if and only if the largest clique of $G(M)$ has cardinality less than $r+2$.  However, determining the cardinality of the largest clique is NP-hard. 
\end{proof}

Besides the connections to graph theory, the monomial case can be used as a tool for reducing the dimensionality of the computations described in Section \ref{sec:equiv}. 

\begin{prop}\label{prop:monocontain}
Let $A \subset S^dV^*$ be any vector space of forms of degree $d,$ and let $M(A)$ denote the span of the monomials that appear as a term with nonzero coefficient in some polynomial in $A$.  Then  $A^{(r)} \subseteq M(A)^{(r)}.$
\end{prop}

\begin{proof}
A monomial differential operator is injective on the set of monomials it does not kill.  Thus, for every monomial $\bfx^{\bal}$ of every polynomial in $A^{(r)}$, and every divisor $\bfx^{\bbe}$ of $\bfx^{\bal}$ there exists a polynomial $f \in A$ such that $\bfx^{\bal - \bbe}$ appears with a nonzero coefficient in $f$.
\end{proof}

Proposition \ref{prop:monocontain} can be useful for computations because the monomial case can be precomputed combinatorially.  Then when applying the algorithms from the previous sections, one can immediately eliminate any polynomials that arise in a partial computation that do not belong to $M(A)^{(r)}$.  Furthermore, the monomial case can be used as a theoretical tool to prove that certain prolongations are, in fact, empty.

\begin{ex}[No 3-way Interaction]
Recall that the toric ideal of the no 3-way interaction model is the kernel of the ring homomorphism
$$\phi_{lmn}:  \kk[x_{ijk}\, \, | \,\, i \in [l], j \in [m], k \in [n] ]  \rightarrow  \kk[a_{ij}, b_{ik}, c_{jk}], \quad  x_{ijk} \mapsto  a_{ij}b_{ik}c_{jk}.$$
The no 3-way interaction model is an example of a log-linear model in statistics.
Giving a complete description of the toric ideals $I_{lmn} = \ker \phi_{lmn}$ is a challenging open problem in algebraic statistics that has been studied by many authors  \cite{Aoki2003, Diaconis1998}.  It is known that the lowest degree of a minimal generator is $4$ and that $A = (I_{lmn})_4$ is spanned by the ${ l \choose 2}{m \choose 2} {n \choose 2}$ binomials that are equivalent to 
$$x_{111}x_{122}x_{212}x_{221} - x_{112}x_{121}x_{211}x_{222}$$
under the natural action of the product of symmetric groups $S_l \times S_m \times S_n$ on indices.

Let $M(A)$ be the space of quartics spanned by the monomials appearing in these binomials.  We will show that $M(A)^{(k)} = 0$ for all $k$.  Proposition \ref{prop:monocontain} then implies that $A^{(k)} = 0 $ for all $k$.  Since the prolongation of a prolongation is a prolongation, it suffices to show that $M(A)^{(1)}= 0$.  This, in turn, is equivalent to showing that there is no monomial of degree five in $\kk[\bfx]$ which is divisible by five distinct monomials from $M(A)$.  However, if we are given any three variables that are part of a monomial in $M(A)$, there is a unique way to complete it to a monomial in $M(A)$, which guarantees that no degree five monomials of the desired type could exist.

Applying Theorem \ref{thm: sec}, we have shown that $A^{(3(r-1))}=0$, and hence that the degree $3r+1$ piece of $I^{\{r\}}_{lmn}$ is zero for all $r, l, m $ and $n$.  This implies that these secant ideals cannot be generated in their lowest possible degree.  \qed
\end{ex}

Another useful property of the monomial point of view is that generation by circuits is preserved when taking prolongations.

\begin{defn}
Let $A \subset S^dV^*$ be a vector space of polynomials and $f \in A.$ The \emph{support} of $f$ is the set of monomials that appear with nonzero coefficient.  The polynomial $f$ is a \emph{circuit} of $A$ if any $g \in A$ with ${\rm supp}(g)  \subseteq {\rm supp}(f)$ is a scalar multiple of $f,$ in other words, $f$ has minimal monomial support.  We say that $A$ is \emph{minimally generated} by its circuits if the set of all of its circuits is a basis for $A.$ 
\end{defn}

\begin{rmk}
Circuits are basic objects in matroid theory that generalize linearly dependent sets.
Note that for polynomials there are two natural definitions of circuits.  One is the definition that we have used, where we consider the set of polynomials as a vector space, and the connection to linear algebra is clear.  Another definition of circuit for an ideal $I$ is that a circuit of $f$ is a polynomial such that the set of variables appearing in $f$ is minimal with respect to inclusion among all nonzero polynomials in $I$.  If $I$ is a prime ideal, this leads to the notion of an algebraic matroid.  This is the definition of circuit which appears, for instance, in Chapter 4 of \cite{Sturmfels1996}, but this is \emph{not} the notion of circuit that we mean. 
\end{rmk}

Note that $A \subset S^dV^*$ is minimally generated by its circuits if and only if the monomial support of any two circuits  with distinct support  are disjoint.  Indeed, suppose that $f$ and $g$ are two circuits that contain the monomial $\bfx^{\bal}$ with coefficient one.  Then the polynomial $f - g \in A$ does not contain $\bfx^{\bal}$.  The new polynomial can be written as a linear combination of circuits with support in ${\rm supp}(f-g)$.  But this implies that either $f$ or $g$ is not needed as a minimal generator.  We show below that the property of being minimally generated by circuits is preserved under prolongation.

\begin{prop}\label{prop:circuits}
If $A$ is minimally generated by its circuits, then so are the prolongations $A^{(r)}$.
\end{prop}

\begin{proof}
It suffices to show that if $\bfx^{\bal}$ is a monomial that is in the support of some circuit of $A^{(r)}$, there is no other circuit of $A^{(r)}$ containing $\bfx^{\bal}$ in its support.  Suppose to the contrary that there are two circuits $f$ and $g$ that contain a monomial in common.  Let $S$ be the set of monomials appearing in both $f$ and $g$.  Let $\bfx^{\bal}$ be any monomial in $S,$ and let $\bfx^{\bbe}$ be any divisor of $\bfx^{\bal}$ of degree $r$.   The derivatives $\frac{\partial^r f}{\partial \bfx^{\bbe}}$ and $\frac{\partial^r g}{\partial \bfx^{\bbe}}$ are thus nonzero and in $A.$ Moreover, they must both contain a multiple of the same circuit $h$ that contains $\bfx^{\bal - \bbe}$, because $A$ is minimally generated by its circuits.  This means that if $\bfx^{\bga}$ appears in $h$, then $\bfx^{\bbe + \bga}$ appears in both $f$ and $g$, and will belong to $S$.  

Now let $f_S$ and $g_S$ be the polynomials obtained from $f$ and $g$ by summing those terms corresponding to elements of $S$.  We will argue that if $\bx^{\bbe}$ has degree $r$, then $\frac{\partial^r f_S}{\partial \bx^{\bbe}}$ is in $A.$  Indeed, by the argument in the preceding paragraph, if $\bx^{\bbe}$ divides an element of $S,$ we get an element of $A,$ and otherwise $\frac{\partial^r f_S}{\partial \bx^{\bbe}} = 0.$  However, since $f$ and $g$ were circuits, this implies that $f_S = f$ and $g_S = g$.  Since $f$ and $g$ are circuits with the same support, they must be nonzero multiples of each other.
\end{proof}

Proposition \ref{prop:circuits} is useful in special cases for proving that we have determined a complete generating set for a particular prolongation.

\begin{ex} 
Let $I$ be the ideal of the Segre embedding of $\pp^{n_1 -1} \times \pp^{n_2 - 1}$ into $\pp^{n_1n_2 -1}$.    The ideal $I$ is generated by the $2 \times 2$ minors of the generic matrix $X = (x_{ij})$.    Let $A$ be the space of quadrics spanned by these $2 \times 2$ minors.  Note that each $2 \times 2$ minor is a circuit, and these circuits have disjoint monomial support.  By Proposition \ref{prop:circuits} $A^{(r)}$ is generated by circuits which have disjoint monomial support.

The monomials appearing in the $2 \times 2$ minors generating $I$ are precisely the monomials $x_{i_1j_1}x_{i_2j_2}$ such that $i_1 \neq i_2$ and $j_1 \neq j_2$.  The cliques in the resulting graph $G(A, r)$  are the monomials of the form $x_{i_1j_1} \cdots x_{i_{r+2}j_{r+2}}$ such that $i_k \neq i_l$ and $j_k \neq j_l$ for all $k \neq l$.  Each such monomial is a term of a unique $(r+2) \times (r+2)$ minor.  Each $(r+2) \times (r + 2)$ minor belongs to $A^{(r)}$ which can be verified by differentiating the Laplace expansion of the determinant.
As each minor is a circuit, we deduce that these $(r+2) \times (r+2)$ minors span the prolongation.  \qed
\end{ex}

%%%%%%%%%%%%%%%%%%%%%%%%%%%%%%%%%%%%%
%%%%%%%%%%%%%%%%%%%%%%%%%%%%%%%%%%%%%
%%%%%%%%%%%%%%%%%%%%%%%%%%%%%%%%%%%%%
%%%%%%%%%%%%%%%%%%%%%%%%%%%%%%%%%%%%%

\section{Prolongations and Secant Varieties}\label{sec:secant}

In this section we will explain the relationship between prolongation and secant varieties.  The proof of Lemma 2.2 in \cite{LM1} goes through in a more general setting, and our proof of Theorem \ref{thm: sec1} follows along these lines.  Although we could also use the ideas of the proof of  Lemma 2.2 in \cite{LM1} to prove Theorem \ref{thm: sec2}, we give an alternate and simpler proof appealing to the computation of joins of ideals.

\begin{thm}\label{thm: sec1}
Suppose that $X \subseteq \pp^{n-1}$ is a variety over an algebraically closed field, and $I = I(X).$  Let $A = I_d$ . Then $A^{((r-1)(d-1))}$ is contained in the ideal of the $r$-th secant variety of $X.$
\end{thm}
\begin{proof}
Suppose that $F \in A^{((r-1)(d-1))}$ so that
\[\deg F = (d-1)(r-1)+d = dr -r-d+1+d = dr-r+1 = r(d-1)+1.\]
A general point on the $r$-th secant variety of $X$ is the span of $r$ points of $X.$  So, let ${\bf v} = t_1 {\bf v}_1+\cdots+t_r{\bf v}_r$ where the $t_i$ and $\bfv_i$ are indeterminates.  We will show that for any specialization of the $\bfv_i$ to points of $X,$ and  $t_i \in \kk,$ $F({\bf v}) = 0.$ 
Since 
\[F({\bf x}) = \frac{1}{(r(d-1)+1)!} {\bf F}({\bf x}, \ldots,{\bf x}),
\]
$F({\bf v}) = 0$ if and only if ${\bf F}({\bf v}, \ldots, {\bf v}) = 0.$

The point now is that the polarization ${\bf F}(\bx_1, \ldots, \bx_{r(d-1)+1})$ is linear in each set of variables $\bx_i.$  This implies that

\[
{\bf F}(\bx_1, \ldots, \bx_{i-1}, {\bf v} , \bx_{i+1}\ldots, \bx_{r(d-1)+1})
= \sum_{j=1}^r t_j{\bf F}(\bx_1, \ldots, \bx_{i-1}, \bfv_j  , \bx_{i+1},\ldots, \bx_{r(d-1)+1}).
\]

Therefore, we see that 

\[{\bf F}(\bfv, \ldots, \bfv)= \sum_{\bbe \in \nn^r, |\bbe| = r(d-1)+1}
{r(d-1)+1 \choose \bbe} \bft^{\bbe}{\bf F}(\bfv_1^{\beta_1}, \ldots, \bfv_{r}^{\beta_{r}}),\]
 where $\bfv_i$ is repeated $\beta_i$ times.  For each $\bbe$ in the sum, $|\bbe| = r(d-1)+1,$ implies that some $\beta_i \geq d.$  Therefore, ${\bf F}(\bfv_1^{\beta_1}, \ldots, \bfv_{r}^{\beta_{r}})$ can be written as a polynomial whose coefficients have degree $d$ in $\bfv_i.$ Since $F \in A^{((r-1)(d-1))},$  every coefficient of a monomial in the ${\bf y}$-variables in the polynomial 
${\bf F}({\bf x}, \ldots, {\bf x}, {\bf y}_1, \ldots, {\bf y}_{(r-1)(d-1)})$ is in $A$ by part (3) of Lemma \ref{lem: polar2}.  Therefore, each of these degree $d$ coefficients is in $A$ (written in the $\bfv_i$-variables).  Thus, if we specialize all of the $\bfv_j$ to points of $X,$ ${\bf F}(\bfv_1^{\beta_1}, \ldots, \bfv_{r}^{\beta_{r}}) = 0.$  We conclude that ${\bf F}(\bfv, \ldots, \bfv) = 0.$
\end{proof}

We also have the partial converse if we know that the ideal of $X$ does not contain any forms of degree $<d.$

\begin{thm}\label{thm: sec2}
Suppose that $X \subset \pp^{n-1}$ is a variety  over an algebraically closed field and that no form of degree $\leq d-1$ vanishes on $X.$  If $m = r(d-1)+1,$  then  $I(X^{\{r\}})_m = A^{((r-1)(d-1))}.$
\end{thm}

To prove the theorem, we collect some general definitions and results about secants and joins of ideals.
Given a collection of ideals $I_1, \ldots, I_r \subseteq \kk[\bfx]$, their \emph{join} is obtained by constructing the new ideal 
$$I_1 \ast \cdots \ast I_r  = \left(  I_1(\bfy_1) + \cdots + I_r(\bfy_r) + \left<  x_j - \sum_{i =1}^r y_{ij} \mid j \in [n] \right>   \right)  \bigcap  \kk[\bfx]$$
where $\bfy_i$ denotes the vector of variables $\bfy_i = (y_{i1}, \ldots, y_{in})$, $I_i(\bfy_i)$ denotes the ideal $I_i$ with variable $y_{ij}$ substituted for variable $x_j$, and the big ideal in parentheses is contained in the ring $\kk[\bfx, \bfy_1, \ldots, \bfy_r]$.  The $r$-fold join of $I$ with itself is the $r$-th secant ideal $I^{\{r\}} = I \ast \cdots \ast I$.  If $I_1, \ldots, I_r$ are saturated homogeneous radical ideals over an algebraically closed field, then $I_1 \ast \cdots \ast I_r$ is the homogenous ideal representing the embedded join of the projective varieties $V(I_1), \ldots, V(I_r)$.  The secant ideal $I^{\{r\}}$ is the vanishing ideal of the $r$-th secant variety of $V(I)$.  Note that the join construction is commutative and associative, and it respects containments.  Proposition \ref{prop:simis}, due to Simis and Ulrich, puts a restriction on the degrees of forms which may appear in the secant ideal, and Proposition \ref{prop:sull} describes how symbolic powers may be computed via the join operation.

\begin{prop}\cite{Simis2000}\label{prop:simis}
If $I \subseteq \left< x_1, \ldots, x_n \right>^d$ then
$$I^{\{r\}}  \subseteq  \left<  x_1, \ldots, x_n \right>^{r(d-1) + 1}.$$
\end{prop}
 
\begin{prop}\cite{Sullivant2006}\label{prop:sull}
Suppose  that $\kk$ is algebraically closed.  If $I$ is a radical ideal, then $I^{(r)} = I \ast \left<x_1 \ldots, x_n \right>^r$.
\end{prop}

These results lead us to a more general result involving containments of secant ideals in symbolic powers.  Theorem \ref{thm: sec2} reflects information about one graded piece of this containment.

\begin{prop}\label{prop:symbpower}
Suppose that $\kk$ is algebraically closed, $I$ is radical, and $I \subseteq \left< x_1, \ldots, x_n \right>^d$.  Then
$$I^{\{r\}}  \subset  I^{((r-1)(d-1) +1)}.$$
\end{prop}

\begin{proof}
This follows by the chain of containments
$$I^{\{r\}}  =  I \ast I^{\{r-1\}}  \subseteq  I  \ast \left< x_1, \ldots, x_n \right>^{(r-1)(d-1) + 1}    = I^{((r-1)(d-1) + 1)}.$$
The first equality is by the associativity of the join, the second containment follows because joins respect containment together with Proposition \ref{prop:simis} and the third equality follows by Proposition \ref{prop:sull}.
\end{proof}

\noindent {\em  Proof of Theorem \ref{thm: sec2}:}
This is a direct consequence of Proposition \ref{prop:symbpower} and Corollary \ref{cor:symbpower}.  \qed

\noindent  {\em  Proof of Theorem \ref{thm: sec}:}
This is a direct consequence of Theorems \ref{thm: sec1} and \ref{thm: sec2}.  \qed

\begin{rmk}
The proofs of Theorems \ref{thm: sec1} and \ref{thm: sec2} can be extended to the nonreduced case.  To do this requires the replacement of the symbolic power with the differential power, and some more algebraic reasoning in the proof of Theorem \ref{thm: sec1}.  We have only included the proof of the reduced case because it is, by far, the most interesting. 
\end{rmk}

%%%%%%%%%%%%%%%%%%%%%%%%%%%%%%%%%%%%%
%%%%%%%%%%%%%%%%%%%%%%%%%%%%%%%%%%%%%
%%%%%%%%%%%%%%%%%%%%%%%%%%%%%%%%%%%%%
%%%%%%%%%%%%%%%%%%%%%%%%%%%%%%%%%%%%%
\section{Application to the Binary Symmetric Model}\label{sec:phylo}

As mentioned in the introduction, one recent motivation for the detailed study of equations vanishing on secant varieties comes from algebraic statistics, where secant varieties correspond to statistical models called mixture models.  Our goal in this section is to illustrate how  prolongations can be used to derive some nontrivial algebraic constraints on mixture models in situations where it seems difficult to prove directly that the same equations belong to the secant ideal.  In particular, we explore this problem for some models arising in phylogenetics.

To give our description of equations in  the prolongation, we first need to describe the space of quadrics which generate the ideal of the phylogenetic model.  The bulk of this description can be found in \cite{Sturmfels2005}.  First we describe the variables of the toric ideal.  Let $T$ be an unrooted trivalent tree (that is, each vertex of the tree that is not a leaf has degree three) with $n$ leaves.  The ideal of the phylogenetic model $I_T$ lives in the polynomial ring in $2^{n-1}$ indeterminates: 
$$ \kk[q] :=  \kk[q_{\bf i} \, \, | \, \, {\bf i} \in (\zz/2\zz)^n_{\rm even} ]$$
where $(\zz/2\zz)^n_{\rm even}$ is the group of binary strings of length $n$ with sum zero in $\zz/2 \zz$.  These $q$ indeterminates are the Fourier transform of the natural probability coordinates (see \cite{Sturmfels2005}).

The ideal $I_T$ is generated by determinantal quadrics.  Specifically, each internal edge $e$ of the tree induces a \emph{split} of the leaves into two disjoint sets, $A|B$.  The indeterminates are also partitioned into two disjoint sets, namely,   the sets  
$$\{ q_{\bfi} \, \, | \, \, \sum_{j \in A} i_j  = 0 \in \zz/2\zz \} \mbox{  and  }  \{ q_{\bfi} \, \, | \, \, \sum_{j \in A} i_j  = 1 \in \zz/2\zz  \}. $$
These two sets of indeterminates fit into two $2^{|A| -1} \times 2^{|B|-1}$ matrices, $M^e_0$ and $M^e_1$, whose rows are indexed by the strings $\bfi_A$ and whose columns are indexed by $\bfi_B$.  The toric ideal $I_T$ is generated by the set of all $2 \times 2$ minors of the matrices $M^e_0$ and $M^e_1$ as $e$ ranges over all the internal edges of $T$.    Let $A_T$ denote the space of quadrics generated by the determinants described above.

\begin{ex}
For instance if $T$ is the trivalent tree with four leaves with unique internal split $12|34$ then we have
$$M^e_{0}  =  \begin{pmatrix}
q_{0000} &  q_{0011} \\
q_{1100} & q_{1111} \end{pmatrix}  \quad \quad
M^e_1 = \begin{pmatrix}
q_{0101} & q_{0110} \\
q_{1001} & q_{1010} \end{pmatrix} $$
and $I_T$ is a complete intersection of quadrics. \qed
\end{ex}
 
While the description given thus far is rich enough to describe generators of the ideals $I_T$, we need a more involved combinatorial description of the indeterminates in the ring $\kk[q]$ to provide a characterization of the polynomials in the prolongation.  The crucial observation is that associated to each indeterminate $q_{\bfi}$ is a labeling of all edges in the tree $T$, by zeros and ones.  An edge gets the label  $\sum_{j \in A(e)} i_j \in \zz/2 \zz$ where $A(e)$ is one part of the split induced by the edge $e$.  Note that this labeling naturally corresponds to a set of disjoint paths through the tree $T$ such that the end points of every path are leaves of the tree.  Conversely, every such set of disjoint paths is the associated labeling of some indeterminate $q_{\bfi}$.  Thus, for each such labeling $L$ we get an indeterminate $q_L$.  In \cite{Buczynska2006}, these unions of paths are called sockets.  We will use these path indeterminates in the remainder of the section.

Now we wish to describe generators of the prolongation of the space of quadrics we have described, which we do in terms of the path indeterminates from above.    A \emph{frame} $F$ is a partial labeling of the tree $T$ where the labels have been assigned to a trivalent subtree $T(F)$ of $T$.  The frame has \emph{active edges} $a(F)$ which are the leaves of $T(F)$ that are not leaves of $T$.  Each active edge $e$ induces a subtree $T_e(F) \subset T$ consisting of all edges on the side of $e$ that does not contain $T(F)$.  Let $L^e$ denote the set of all possible labelings of $T_e(F)$ that are compatible with the frame $F$ (that is, can be completed to a variable).  If $F$ is a frame and $e$ is an edge of $F$ that has been assigned, let $F(e)$ be the label assigned to the edge $e$.

\begin{defn}
A collection of frames $F_1, \ldots, F_d$ together with a function $e(\cdot, \cdot) : {[d] \choose 2} \rightarrow E(T)$ is \emph{compatible} if $e(\cdot , \cdot)$ satisfies:
\begin{enumerate}
\item  For all $(i,j) \in {[d] \choose 2}$, 
$ e(i,j)  \in a(F_i) \cap a(F_j) \mbox{ with }  F_i(e(i,j)) = F_j(e(i,j))$;
\item  If $e(i,j) = e(j,k)$ then $e(i,j) = e(i,k)$; and
\item  For all $j \in [d]$,  $ \cup_{i \neq j}  e(i,j)    =  a(F_j)$.
\end{enumerate}
\end{defn}

The function $e(\cdot, \cdot)$ determines an equivalence relation on the set of pairs $(F_i, e)$ with $e \in a(F_i)$ where $(F_i, e)$ is defined to be equivalent to itself, and if $i\neq j,$ then we define $(F_i, e_1) \sim (F_j, e_2)$ if $e(i,j) = e_1 = e_2$.  Let $E$ denote such an equivalence class and let $C(E) \subset L^e $ be a set of $|E|$ distinct labelings of $T_e(F)$ compatible with the $F_i \in E$.  Given all these data (in particular, the frames and the labeling sets $C(E)$) we define a polynomial of degree $d$.  To do this, fix a particular base ordering on each of the sets $C(E)$.  Now for each $E$, take some permutation of $C(E)$.  This set of permutations can be used to complete all the frames $F_1, \ldots, F_d$ in a unique ordered way.  This is accomplished by adding the first element of each $C(E)$ to the frame $F_j$ that appears first in the equivalence class $E$ and so on.  Thus each set of permutations yields a monomial in the $q_L$.  Denote by $P(F_1, \ldots, F_d; C(E_1), \ldots, C(E_k))$ the signed sum of all such monomials where the sign of a monomial is the product of the signs of the each permutations used to form the monomial.

\begin{ex}[The six-leaf snowflake]\label{ex:snow}
Let $T$ be the six-leaf tree depicted below and let $A$ be the span of the quadratic binomials which generate $I_T.$

\[
\includegraphics[width = 4cm]{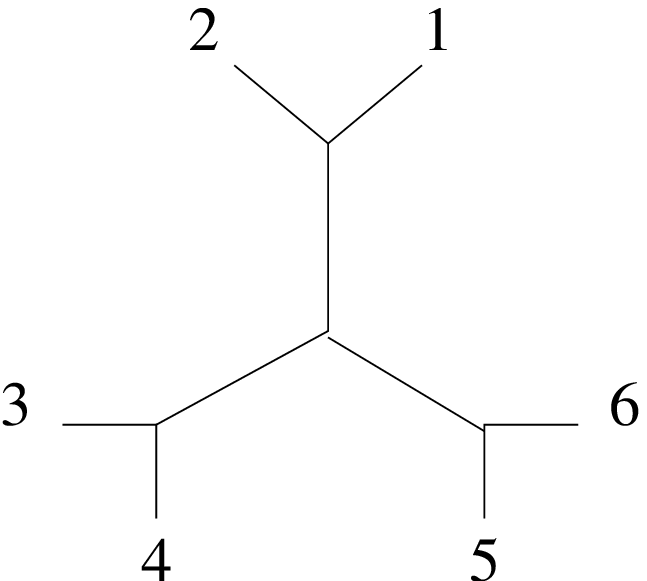}
\]

Computing $A^{(1)}$ and $A^{(2)}$ with \emph{Macaulay 2} \cite{GS}, we see that $A^{(1)}$ is spanned by 32 8-term cubics and $A^{(2)}$ is spanned by a single 64-term quartic form.  As we will see, in both cases, the construction described above yields the entire prolongation.  

To construct a cubic, we need to choose three frames with compatible labelings.  In the diagram below our choice of frame is depicted in bold. 

\[
\includegraphics{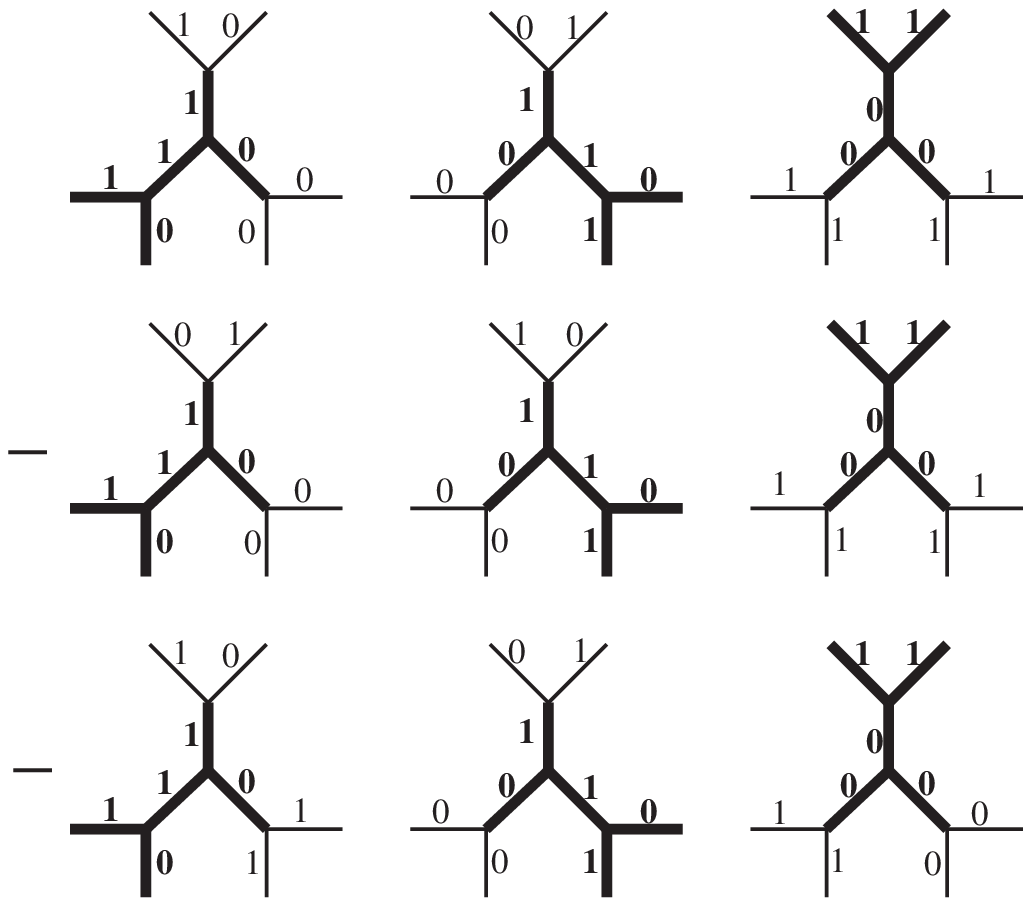} 
\]

 We show triples of trees that correspond to the first three terms of the 8-term cubic

\begin{align*}
  q_{011000} q_{100010} q_{111111}  & - q_{101000}q_{010010}q_{111111} \\
-q_{011011} q_{100010} q_{111100}  & + q_{101011}q_{010010} q_{111100} \\
-q_{011000}q_{101110} q_{110011}  &  + q_{101000}q_{011110}q_{110011} \\
+q_{011011}q_{101110}q_{110000} & - q_{101011}q_{011110}q_{110000}
\end{align*}

We get 32 such cubics because there are 4 ways of choosing 3 distinct inner edges and two ways of labeling each of the pairs of cherries attached to the three inactive edges on the frames.

Our 64-term quartic is constructed by choosing the 4-tuple consisting of our 4 distinct choices of frames and completing each edge-labelling in any way allowed. There is only one way to define the function $e(\cdot, \cdot)$ and each edge on each frame is active.
\end{ex}

\begin{ex}[The six-leaf caterpillar]\label{ex:cat}
Let $T$ be the six-leaf caterpillar-shaped graph below.

\[
\includegraphics[width = 4cm]{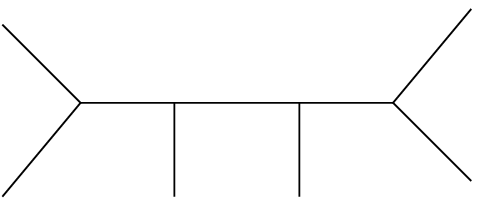}
\]

The corresponding toric variety has ideal $I_T$ generated by the $2 \times 2$ minors of four $2 \times 8$ and two $4 \times 4$ matrices.  Using \emph{Macaulay 2} \cite{GS}, one can see that $A^{(1)}$ is spanned by 32 6-term cubics and that $A^{(2)}$ is spanned by two 24-term quartics.  These forms may be constructed using the methods above although it is easy to see that they are the $3 \times 3$ and $4 \times 4$ minors of the $4 \times 4$ matrices used to define $I_T.$
\end{ex}

\begin{thm}
For any set of compatible frames $\mathcal{ F} = F_1, \ldots, F_d$, compatibility function $e(\cdot, \cdot)$, and completions $\mathcal{ C} = C(E_1), \ldots, C(E_k)$, the polynomial $P(\mathcal{ F} ; \mathcal{ C} )$ is in the prolongation $A_T^{(d-2)}$.
\end{thm}

\begin{proof}
The proof is by induction on $d$.  When $d = 2$,  $P(\mathcal{ F} ; \mathcal{ C} )$ yields the description of a $2 \times 2$ minor of a matrix $M^{e(1,2)}_{F(e(1,2))}$ associated to the unique common active edge in the two frames.  So suppose the statement is true for $d-1$.  The result will follow if we show that the derivative of the degree $d$ polynomial $P(\mathcal{ F}; \mathcal{ C})$ with respect to any variable is the sum of polynomials of the form $P(\mathcal{ F}'; \mathcal{ C}')$ of degree $d-1.$

Let $\mathcal{ F} = F_1, \ldots, F_d$ and $q_L$ be any variable appearing in a monomial in $P(\mathcal{ F}; \mathcal{ C})$.   By our construction of $P(\mathcal{ F}; \mathcal{ C}),$ each occurrence of $q_L$ is associated to a frame in $\mathcal{F}.$   Without loss of generality, assume that $q_L$ arises by completing the labeling of the frame $F_d$.  (It may be associated to other frames as well.)   In each monomial in which $q_L$ appears by completing a labeling of $F_d,$ the $d-1$ other factors come from the frames $F_1, \ldots F_{d-1}$.  Now construct the new set $C(E_1)', \ldots, C(E_k)'$ by removing the elements from $C(E_1), \ldots, C(E_k)$ that are used to make $q_L$.  If any of the sets $C(E_i)'$ are singletons, we can take this single element and modify the appropriate frame $F_l$ to get a new frame, and remove the set $C(E_i)'$ from our list of completions.  Carrying out this procedure yields a set of frames $\mathcal{ F}' = F_1', \ldots, F_{d-1}'$ and a set of completions $\mathcal{ C}' = C(E_1'), \ldots C(E_r')$, such that upon dividing all monomials in $P(\mathcal{ F}; \mathcal{ C})$ that have this particular realization of $q_L$ (associated to the frame $F_d$) we get the polynomial $P(\mathcal{ F}'; \mathcal{ C}')$.  Applying this argument to all realizations of $q_L$ by different frames, we deduce that the derivative of $P(\mathcal{ F}; \mathcal{ C})$ with respect to $q_L$ is the (signed) sum of polynomials $P(\mathcal{ F}'; \mathcal{ C}')$.
\end{proof}

\begin{rmk}
Note that the argument in the preceding proof holds even when $q_L$ appears to a power $>1$, because the coefficient of the derivative will account for the different frames that yield $q_L$.  It is interesting to note that this exceptional case cannot occur, however, because $A$ is generated by polynomials with all squarefree terms.  Thus, $A^{(d-2)}$ is also generated by polynomials with all squarefree terms.
\end{rmk}

%%%%%%%%%%%%%%%%%%%%%%%%%%%%%%%%%%%%%
%%%%%%%%%%%%%%%%%%%%%%%%%%%%%%%%%%%%%
%%%%%%%%%%%%%%%%%%%%%%%%%%%%%%%%%%%%%
%%%%%%%%%%%%%%%%%%%%%%%%%%%%%%%%%%%%%


\begin{thebibliography}{99}

\bibitem{Aoki2003} S.~Aoki and A.~Takemura.   Minimal basis for a connected Markov chain over $3\times 3\times K$ contingency tables with fixed two-dimensional marginals. \emph{Aust. N. Z. J. Stat.} {\bf 45} (2003), no. 2, 229--249.

\bibitem{ACGH} E.~Arbarello, M.~Cornalba, P.~A.~Griffiths and J.~Harris.  \emph{Geometry of Algebraic Curves, Volume I}, Springer-Verlag, New York, 1985. 

\bibitem{magma} W.~Bosma, J.~Cannon, and C.~Playoust. The Magma algebra system.~I.~The user language. {\em J. Symbolic Computation}, {\bf 24} (1997) no.~3--4,  235--25


\bibitem{BCGGG} R.~L.~Bryant, S.~S.~Chern, R.~B.~Gardner, H.~L.~Goldschmidt, P.~A.~Griffiths. \emph{Exterior differential systems}, Springer-Verlag, New York, 1991.

\bibitem{Buczynska2006}  W.~Buczynska and J.~Wisniewski.  On phylogenetic trees--  a geometer's view.  Preprint, 2006.  {\tt math.AG/0601357}


\bibitem{cj} M.~Catalano-Johnson. The homogeneous ideals of higher secant varieties. \emph{J. Pure Appl. Algebra} {\bf 158} (2001), no.~2--3, 123--129. 

\bibitem{Diaconis1998}  P.~Diaconis and B.~Sturmfels.  Algebraic algorithms for sampling from conditional distributions.  \emph{Annals of Statistics} {\bf 26} (1998), no. 1, 363--397.

\bibitem{E}  D.~Eisenbud, \emph{Commutative Algebra: with a View Toward Algebraic Geometry}, Springer-Verlag, New York, 1995.

\bibitem{FH} W.~Fulton and J. Harris \emph{Representation Theory}, Springer-Verlag, New York, 1991.

\bibitem{G}  P. Griffiths. Some aspects of exterior differential systems. {\em Proceedings of Symposia in Pure Mathematics}. {\bf 53} (1991), 151--173.

\bibitem{GS} D.~Grayson and M.~Stillman, \emph{Macaulay 2, a software system for research in algebraic geometry}, available at {\tt http://www.math.uiuc.edu/Macaulay2/}


\bibitem{IL}  T.~Ivey and J.M.~Landsberg, \emph{Cartan for Beginners: Differential Geometry Via Moving Frames and Exterior Differential Systems}, Graduate Studies in Math. {\bf 61}, American Mathematical Society, Providence, 2003.


\bibitem{LM1} J.M. Landsberg and L. Manivel, On the projective geometry of rational homogeneous varieties. {\em Comment. Math. Helv.}, {\bf 78} (2003), no.1, 65--100.

\bibitem{LM2} J.M. Landsberg and L. Manivel.  On the ideals of secant varieties of Segre varieties, {\em Found. Comput. Math.} (2004), no.~4, 397--422. 

\bibitem{Simis2000}  A.~Simis and B.~Ulrich. On the ideal of an embedded join. \emph{Journal of Algebra}  {\bf 226}  (2000), no.~1, 1--14.
\bibitem{Sturmfels1996}  B.~Sturmfels.  \emph{Gr\"obner Bases and Convex Polytopes}. University Lecture Series {\bf 8},  American Mathematical Society, Providence, RI, 1996.

\bibitem{Sturmfels2005}  B.~Sturmfels and S.~Sullivant.  Toric ideals of phylogenetic invariants.  {\em Journal of Computational Biology}  {\bf 12}  (2005) 204-228.

\bibitem{Sturmfels2006}  B.~Sturmfels and S.~Sullivant.  Combinatorial secant varieties. \emph{ Quarterly Journal of Pure and Applied Mathematics} {\bf 2} (2006), no.~3, 285-309.  

\bibitem{Sullivant2006}  S.~Sullivant.  Combinatorial symbolic powers.  Preprint, 2006.  {\tt math.AC/0608542}

\bibitem{v} P.~Vermeire, Some results on secant varieties leading to a geometric flip construction. \emph{Compositio Math.} {\bf 125} (2001), no.~3, 263--282.

\bibitem{W} J.~Weyman, \emph{Cohomology of Vector Bundles and Syzygies}, Cambridge University Press,  New York, 2003.

\bibitem{Weyl} H.~Weyl, \emph{Classical Groups: Their Invariants and Representations}, Princeton University Press, Princeton, 1939.
\end{thebibliography}
\end{document}